\documentclass{article}
\usepackage{latexsym}
\usepackage{amstext}
\usepackage{amsmath}
\usepackage{amssymb}
\usepackage{amsthm}

\newtheorem{theorem}{Theorem}
\newtheorem{lemma}[theorem]{Lemma}

\newtheorem{question}[theorem]{Question}
\theoremstyle{definition}
\newtheorem{definition}[theorem]{Definition}

\def \co {\mathcal{O}}

\def \kk {\bar{k}}
\def \pn {\mathbb{P}^n}
\DeclareMathOperator{\Tr}{Tr}
\DeclareMathOperator{\Span}{Span}
\DeclareMathOperator{\Gal}{Gal}

\begin{document}
\bibliographystyle{amsplain}
\title{On the Zariski-Density of Integral Points on a Complement of Hyperplanes in $\mathbb{P}^n$}
\author{Aaron Levin\\adlevin@math.brown.edu}
\date{}
\maketitle
\begin{abstract}
We study the $S$-integral points on the complement of a union of hyperplanes in projective space, where $S$ is a finite set of places of a number field $k$.  In the classical case where $S$ consists of the set of archimedean places of $k$, we completely characterize, in terms of the hyperplanes and the field $k$, when the ($S$-)integral points are not Zariski-dense.
\end{abstract}
\section{Introduction}
Let $k$ be a number field and $S$ a finite set of places of $k$ containing the archimedean places.  Let $Z$ be a closed subset of $\pn$, defined over $k$, that is a finite union of hyperplanes over $\kk$.  We study the problem of determining when there exists a Zariski-dense set $R$ of $S$-integral points on $\mathbb{P}^n\backslash Z$.   We give a complete answer to this problem when $\co_{k,S}=\co_k$ is the usual ring of integers, i.e., when $S=S_\infty$ consists of the set of archimedean places of $k$.  For arbitrary $S$ the problem does not appear to have a simple answer, but in the last section we discuss some partial results and reformulations of the problem.

The related problem of determining when $R$ must be a finite set was solved by Evertse and Gy\"ory \cite{Ev} once $k$ is sufficiently large (e.g., the hyperplanes are all defined over $k$).  In the connected topic of solutions to norm form equations, Schmidt \cite{Sc2}\cite{Sc} has given necessary and sufficient conditions for finiteness.  The general problem of determining the possible dimensions of $R$, for any $k$, $S$, and $Z$, seems to be difficult.

\section{Definitions}
Let $k$ be a number field and $S$ a finite set of places of $k$ containing the archimedean places.  Let $\co_{k,S}$ denote the ring of $S$-integers of $k$.
\begin{definition}
If $Z$ is a subset of $\pn$ defined over $k$, we call a set $R\subset \pn\backslash Z(k)$ a {\it set of $S$-integral points on $\pn \backslash Z$} if for every regular function $f$ on $\pn \backslash Z$ defined over $k$ there exists $a\in k^*$ such that $af(P)\in \co_{k,S}$ for all $P\in R$.
\end{definition}
Equivalently, if $Z$ is a hypersurface, $R$ is a set of $S$-integral points on $\pn \backslash Z$ if there exists an affine embedding $\pn \backslash Z\subset \mathbb{A}^N$ such that every $P\in R$ has $S$-integral coordinates.

Recall that an archimedean place $v$ of $k$ corresponds to an embedding of $k$ into the complex numbers $\sigma:k\to \mathbb{C}$.  We define $v$ to be real if $\sigma(k)\subset \mathbb{R}$ and define $v$ to be complex otherwise.  With this terminology we can define the following types of fields.
\begin{definition}
Let $k$ be a number field.  Then
\begin{enumerate}
\item  We call $k$ a totally real field if all of its archimedean places are real.
\item  We call $k$ a totally imaginary field if all of its archimedean places are complex.
\item  We call $k$ a complex multiplication (CM) field if it is a totally imaginary field that is a quadratic extension of a totally real field.
\item  We say that an extension $M$ of $k$ contains a CM subfield over $k$ if there exists a CM field $L$ with maximal real subfield $L'$ (over $\mathbb{Q}$) such that $k\subset L' \subset L \subset M$.
\end{enumerate}
\end{definition}
Note that in our terminology, if $M$ is a CM field then $M$ does not contain a CM subfield over itself because of the condition on the maximal real subfield.
\section{Main Theorem}
Our main theorem gives a complete characterization of when there exists a Zariski-dense set of ($S_\infty$-)integral points on a complement of hyperplanes.
\begin{theorem}
\label{Hyper}
Let $Z\subset \pn$ be a closed subset defined over $k$ that is a geometric finite union of hyperplanes, i.e., $Z=\cup_{i=1}^m H_i$ over $\kk$ where the $H_i$ are distinct hyperplanes defined over $\kk$.  Let $L_i$ be a linear form defining $H_i$ over its minimal field of definition $M_i$.  Let $S=S_\infty$, the set of archimedean places of $k$.  Then there does not exist a Zariski-dense set of $S$-integral points on $\pn \backslash Z$ if and only if one of the following conditions holds:
\begin{enumerate}
\item  The linear forms $L_1,\ldots,L_m$ are linearly dependent.
\item  $\co_{k,S}^*=\co_k^*$ is finite and $Z$ has more than one irreducible component over $k$.
\item  Some $M_i$ contains a CM subfield over $k$.
\end{enumerate}
\end{theorem}
\begin{proof}
We prove the ``if" direction first.  Suppose that (a) holds.    Without loss of generality, we can extend $k$ so that each $L_i$ is defined over $k$.  It suffices to prove our assertion in the case that $\{L_1,\ldots,L_m\}$ is a minimal linearly dependent set, that is no proper subset is linearly dependent.  In that case $\sum_{i=1}^{m-1} c_iL_i=c_mL_m$ for some choice of $c_i\in k^*, i=1,\ldots,m$.  Let $R$ be a set of $S$-integral points on $\pn \backslash Z$.  If $i \in \{1,\ldots, m\}$, then all of the poles of $L_i/L_m$ lie in $Z$ and so there exists an $a\in k^*$ such that $af$ takes on integral values on $R$.  Since the poles of $L_m/L_i$ also lie in $Z$, the same reasoning applies to $L_m/L_i$.  Therefore $L_i/L_m(R)$ is contained in the union of finitely many cosets of the group of units $\co_{k,S}^*$.  By enlarging $S$ we can assume without loss of generality that $\frac{c_iL_i}{c_mL_m}(P)$ is a unit for all $P\in R$ and $i=1,\ldots,m$.  We now apply the $S$-unit lemma.
\begin{lemma}[$S$-unit Lemma]
Let $k$ be a number field and $n$ a positive integer.  Let $\Gamma$ be a finitely generated subgroup of $k^*$.  Then all but finitely many solutions of the equation
\begin{equation*}
u_0+u_1+\cdots+u_n=1, u_i\in \Gamma
\end{equation*}
satisfy an equation of the form $\sum_{i\in I}u_i=0$, where $I$ is a subset of $\{0,\ldots,n\}$.
\end{lemma}
Since $\sum_{i=1}^{m-1} \frac{c_iL_i}{c_mL_m}(P)=1$ for all $P\in R$,  by the $S$-unit lemma it follows that each $P\in R$ either belongs to one of the hyperplanes defined by $\sum_{i\in I}c_iL_i=0$ for some subset $I\subset \{1,\ldots,m-1\}$ (this equation is nontrivial by the minimality of the linear dependence relation) or it belongs to a hyperplane defined by $c_iL_i=tc_mL_m$, for some $t\in T$, where $T \subset \co_k^*$ is a finite subset containing the elements that appear in the finite number of exceptional solutions to the $S$-unit equation $\sum_{i=1}^{m-1} x_i=1$.  In particular, $R$ is not Zariski-dense.

Suppose that (b) holds.  Let $R$ be a set of $S$-integral points on $\pn \backslash Z$.  Let $Z_1$ and $Z_2$ be two distinct irreducible components of $Z$ defined over $k$, respectively, by homogeneous polynomials $f$ and $g$.  Let $h=f^{\deg g}/g^{\deg f}$.   Since both $h$ and $1/h$ are regular on $\pn\backslash Z$, by our earlier argument $h(R)$ is contained in the union of finitely many cosets of $\co_{k}^*$.  By our assumption on $\co_k^*$, $h(R)$ is a finite set.  This implies that $R$ is contained in the union of finitely many hypersurfaces of the form $f^{\deg g}=ag^{\deg f}, a\in k$, and so $R$ is not Zariski-dense.

Suppose that (c) holds.  It suffices to prove our assertion in the case that $Z$ is irreducible over $k$, $H_1$ has minimal field of definition $M$, and $M$ contains a CM subfield $L$ over $k$.  From what we have already proven, we can assume that the linear forms defining the hyperplanes are linearly independent.  It follows from the fact that $Z$ is irreducible that $[M:k]=m$.  Let $\alpha_0,\ldots,\alpha_{m-1}\in \co_M$ be a basis for $M$ over $k$.  Under our assumptions, after a $k$-linear change of variables (a projective $k$-automorphism of $\mathbb{P}^{n}$), we can take $Z$ to be defined by $N^M_k(x_0\alpha_0+\cdots+x_{m-1}\alpha_{m-1})=0$, where $N^M_k$ is the norm from $M$ to $k$, and the embeddings of $M$ act on each $x_i$ trivially.  From the defining equation for $Z$, it suffices to prove the case $n=m-1$.
\begin{lemma}
\label{units}
Let $Z\subset \mathbb{P}^{m-1}$ be defined by $N^M_k(x_0\alpha_0+\cdots+x_{m-1}\alpha_{m-1})=0$.  Let $R$ be a set of $S$-integral points on $\mathbb{P}^{m-1}\backslash Z$.  Let $S_M$ be the set of places of $M$ lying above places of $S$.  There exist a finite number of elements $\beta_1,\ldots,\beta_r\in M$ such that every $P\in R$ has a representative $P=(x_0,\ldots,x_{m-1})\in \mathbb{P}^{m-1}(\co_{k,S})$ with $\sum_{i=0}^{m-1}x_i\alpha_i\in \beta_j\co_{M,S_M}^*$ for some $j$.
\end{lemma}
\begin{proof}
By the definition of $R$ being an $S$-integral set of points on $\mathbb{P}^{m-1}\backslash Z$, for any monomial $p$ in $x_0,\ldots,x_{m-1}$ of degree $m$, there exists a constant $c_p\in k^*$ such that $c_p p/N^M_k(\sum_{i=0}^{m-1}x_i\alpha_i)$ takes on $S$-integral values on $R$.  Therefore there exists a constant $C\in k^*$ such that for all points $(x_0,\ldots,x_{m-1})\in R$,
\begin{equation}
\label{HyperGE}
\left(N^M_k\left(\sum_{i=0}^{m-1}x_i\alpha_i\right)\right)|C(x_0,\ldots,x_{m-1})^m
\end{equation}
as fractional ideals of $\co_{k,S}$.  Since the class group of $k$ is finite, there exists a finite set of integral ideals $\mathfrak{A}$ such that for any point of $R$ we can write 
\begin{equation*}
(x_0,\ldots,x_{m-1})=(\beta)\mathfrak{a},
\end{equation*}
where $\beta\in k$ and $\mathfrak{a}\in \mathfrak{A}$.  So dividing (\ref{HyperGE}) by $\beta^m$ on both sides, we see that every point of $R$ has a representative $(x_0,\ldots,x_{m-1})\in \mathbb{P}^{m-1}(\co_{k,S})$ such that (as $\co_{k,S}$ ideals) 
\begin{equation*}
\left(N^M_k\left(\sum_{i=0}^{m-1}x_i\alpha_i\right)\right)|\mathfrak{b}
\end{equation*}
where $\mathfrak{b}$ is some fixed ideal of $\co_{k,S}$ independent of $x_0,\ldots,x_{n-1}$.  Modulo $S'$-units, there are only finitely many solutions $x=\beta_1,\ldots,\beta_r$ to 
\begin{equation*}
\left(N^M_k(x)\right)|\mathfrak{b}, \quad x \in \co_{M,S_M}.
\end{equation*}
The claim then follows.
\end{proof}
We make the following convenient definition.
\begin{definition}
Let $M$ be a finite extension of a field $k$, $[M:k]=n$.  Let $R$ be a subset of $M^*$.  Let $\alpha_0,\ldots,\alpha_{n-1}$ be a basis for $M$ over $k$.  We define $R$ to be a dense subset of $M$ over $k$ if the set $\{(x_0,\ldots,x_{n-1})\in \mathbb{P}^{n-1}(k):\sum_{j=0}^{n-1}x_j \alpha_{j} \in R\}$ is a Zariski-dense subset of $\mathbb{P}^{n-1}$.
\end{definition}
This definition is clearly independent of the basis $\alpha_0,\ldots,\alpha_{n-1}$ that is chosen.  If $R\subset M^*$ is not a dense subset of $M$ over $k$ it is clear that $\alpha R$ for $\alpha \in M^*$ is also not a dense subset of $M$ over $k$, since the corresponding subsets of $\mathbb{P}^n$ differ by a projective automorphism.  Therefore, using Lemma \ref{units}, to finish our claim assuming (c) we need to show that $\co_{M,S_M}^*$ is not a dense subset of $M$ over $k$, where $M$ contains a CM subfield $L$ over $k$.  Let $L'$ be the maximal real subfield of $L$.  Since the totally imaginary field $L$ is a quadratic extension of the totally real field $L'$, by the Dirichlet unit theorem the unit groups of $\co_L$ and $\co_{L'}$ have the same free rank.  It follows that there exists a positive integer $m$ such that if $u\in \co_L^*$ then $u^m\in \co_{L'}^*$.  Let $[L:k]=2l$.  Let $\beta_0,\ldots,\beta_{2l-1}$ be a basis for $L$ over $k$ where $\beta_0,\ldots,\beta_{l-1}$ are real (and are therefore a basis for $L'$ over $k$).  Let $[N^{M}_L(\sum_{i=0}^{n-1}x_i \alpha_i)]^m=\sum_{j=0}^{2l-1}f_i \beta_i$ where the $f_i$ are homogeneous polynomials in $x_0,\ldots,x_{n-1}$.  Since for any $u\in \co_M^*$, $N^M_Lu\in \co_L^*$, we obtain $(N^M_Lu)^m\in \co_{L'}^*$.  Therefore the nontrivial polynomials $f_i$ for $l\leq i \leq 2l-1$ vanish on the set associated to $\co_M^*$ in this basis.  So $\co_M^*$ is not a dense subset of $M$ over $k$.

To prove the other direction of the theorem, suppose that (a), (b), and (c) are all not satisfied.  Let $Z_1,\ldots,Z_r$ be the irreducible components of $Z$ over $k$.  For each $i$, let $M_i$ be the minimal field of definition of some hyperplane in $Z_i$.  Let $d_i=[M_i:k]$ and let $s(i)=\sum_{j=1}^{i-1}d_j$.  Let $\alpha_{0,i},\ldots,\alpha_{d_i-1,i}$ be a basis for $M_i$ over $k$.  Since the linear forms $L_i$ are linearly independent, it follows that after a $k$-linear change of coordinates, $Z$ can be defined by 
\begin{equation*}
\prod_{i=1}^r N^{M_i}_k \left(\sum_{j=0}^{d_i-1}x_{s(i)+j} \alpha_{j,i}\right)=0.
\end{equation*}
Additionally, by assumption, $r=1$ if $\co_{k}^*$ is finite.  We claim that the set
\begin{equation*}
R=\left\{(x_0,\ldots,x_n)\in \mathbb{P}^n(k):\forall i,\sum_{j=0}^{d_i-1}x_{s(i)+j} \alpha_{j,i} \in \co_{M_i}^*, \forall l\geq s(r),x_l\in \co_k\right\}
\end{equation*}
is a dense set of $S$-integral points on $\mathbb{P}^n\backslash Z$.  That $R$ is a set of $S$-integral points on $\mathbb{P}^n\backslash Z$ is clear from our defining equation for $Z$, the fact that norms of units are units, and that if $(x_0,\ldots,x_n)\in R$ then $x_i\in \frac{1}{N}\co_k$ for some fixed $N\in \co_k$.  

When $\co_{k}^*$ is infinite, we first give an argument to reduce our claim to the case $r=1$, where $Z$ is irreducible over $k$.  Suppose that $R$ is not Zariski-dense.  Let $P$ be a nonzero homogeneous polynomial with a minimal number of terms vanishing on $R$.  We also choose such a $P$ with minimal degree.  Since for $(x_0,\ldots,x_n)\in R$, $x_l$ for $l\geq s(r)$ can be chosen in the infinite set $\co_k$ independently of the other $x_i$, it is clear that $P$ does not contain any of the variables $x_l$, $l\geq s(r)$.  After reindexing, we can assume that $x_0$ appears in $P$.  It follows from our minimality assumptions about $P$ and the structure of $R$ that one can specialize the variables $x_{d_1},\ldots,x_n$ to obtain a nonzero polynomial $P'(x_0,\ldots,x_{d_1-1})$ (not necessarily homogeneous) that vanishes on the set 
\begin{equation*}
R'=\left\{(x_0,\ldots,x_{d_1-1})\in \mathbb{A}^{d_1}(k):\sum_{j=0}^{d_1-1}x_{j} \alpha_{j,1} \in \co_{M_1}^*\right\}.  
\end{equation*}
Write $P'=\sum_{i=0}^q P_i'$, where each $P_i'$ is homogeneous of degree $i$.  If $u\in\co_{k}^*$ then $P'_u=P'(ux_0,\ldots,ux_{d_1-1})=\sum_{i=0}^q u^iP_i'$ gives another polynomial that vanishes on $R'$.  Since $\co_{k}^*$ is infinite, we can choose $q+1$ distinct units $u_1,\ldots,u_{q+1}$ of $\co_{k}^*$, and by the invertibility of a Vandermonde matrix, we see that for each $i$, $P_i'\in \Span\{P_{u_1},\ldots,P_{u_{q+1}}\}$.  Therefore if $R$ is not Zariski-dense, we obtain a nonzero homogeneous polynomial that vanishes on $R'$.  Showing that such a homogeneous polynomial does not exist is equivalent to the $r=1$ case of our original claim.  In other words, we have reduced the problem, whether or not $\co_{k}^*$ is finite, to showing that if $M$ does not contain a CM subfield over $k$, $[M:k]=n$, the set
\begin{equation}
\label{Req}
R=\left\{(x_0,\ldots,x_{n-1})\in \mathbb{P}^{n-1}(k):\sum_{j=0}^{n-1}x_j \alpha_{j} \in \co_{M}^*\right\}
\end{equation}
is Zariski-dense, where $\alpha_0,\ldots,\alpha_{n-1}$ is a basis for $M$ over $k$.  In our terminology,  we need to show that $\co_M^*$ is a dense subset of $M$ over $k$.
\begin{theorem}
\label{Hyper2}
Let $M$ be a finite extension of a number field $k$.  The set of units $\co_M^*$ of $\co_M$ is a dense subset of $M$ over $k$ if and only if $M$ does not contain a CM subfield over $k$.
\end{theorem}
We will need the following lemma.
\begin{lemma}
\label{HyperL}
Let $M$ be a finite extension of a number field $k$, $[M:k]=n$.  Let $\sigma_1,\ldots,\sigma_n$ be the embeddings of $M$ into $\mathbb{C}$ fixing $k$.  Let $G$ be a multiplicative subgroup of $M^*$.  Then $G$ is not a dense subset of $M$ over $k$ if and only if there exist nonidentical sequences of nonnegative integers $a_1,\ldots,a_n$ and $b_1,\ldots,b_n$ with $\sum_{i=1}^na_i=\sum_{i=1}^nb_i$ such that 
\begin{equation}
\label{eq}
\prod_{i=1}^n \sigma_i(x)^{a_i}=\prod_{i=1}^n \sigma_i(x)^{b_i}
\end{equation}
 for all $x\in G$.
\end{lemma}
\begin{proof}
Let $\alpha_0,\ldots,\alpha_{n-1}$ be a basis for $M$ over $k$.  Let $R$ be the subset of $\mathbb{P}^{n-1}$ associated to $G$ in this basis.  Suppose that there exist nonidentical sequences of nonnegative integers $a_1,\ldots,a_n$ and $b_1,\ldots,b_n$ with $\sum_{i=1}^na_i=\sum_{i=1}^nb_i$ such that $\prod_{i=1}^n \sigma_i(x)^{a_i}=\prod_{i=1}^n \sigma_i(x)^{b_i}$ for all $x\in G$.  Substituting $x=\sum_{i=0}^{n-1}x_i\alpha_i$ into this equation gives a homogeneous polynomial that vanishes on $R$.  It remains to show that this polynomial is nonzero, or equivalently, that for some $x\in M^*$, $\prod_{i=1}^n \sigma_i(x)^{a_i}\not=\prod_{i=1}^n \sigma_i(x)^{b_i}$.  To see this, we can take for example an $x\in \co_M$ such that $(x)=\mathfrak{p}^q$ for some $q$, where $\mathfrak{p}$ lies above a prime of $k$ that splits completely in $\tilde{M}$, the Galois closure of $M$ over $k$.  Looking at the prime ideal factorization (in $\co_{\tilde{M}}$) of both sides shows that they are unequal.  Therefore $G$ is not a dense subset of $M$ over $k$.

Suppose now that there exists a nonzero homogeneous polynomial vanishing on $R$.  If $x\in G$ and $x=\sum_{i=0}^{n-1}x_i\alpha_i$, $x_i\in k$, then it follows from the fact that $\Tr^M_k(xy)$ is a nondegenerate bilinear form over $k$ that each $x_i$ is a linear form, independent of $x$, in $\sigma_1(x),\ldots, \sigma_n(x)$.  Thus, any nonzero homogeneous polynomial vanishing on $R$ gives rise to a nonzero homogeneous polynomial $P(x_0,\ldots,x_{n-1})$ such that $P(\sigma_1x,\ldots,\sigma_nx)=0$ for all $x\in G$.  Let $P$ be such a polynomial with a minimal number of terms.  Let $c_1\prod_{i=1}^n \sigma_i(x)^{a_i}=c_1\phi_1(x)$ and $c_2\prod_{i=1}^n \sigma_i(x)^{b_i}=c_2\phi_2(x)$ be two monomials appearing in $P(\sigma_1x,\ldots,\sigma_nx)$.  Note that $\sum_{i=1}^na_i=\sum_{i=1}^nb_i$.  Suppose that there exists an $a\in G$ such that $\phi_1(a)\not=\phi_2(a)$.  Let $Q=P(\sigma_1(a)x_0,\ldots,\sigma_n(a)x_{n-1})$.  Since $\phi_1(a)\not=\phi_2(a)$, $Q$ is not a multiple of $P$.  Since $G$ is a group, we have 
\begin{equation*}
P(\sigma_1(a)\sigma_1(x),\ldots,\sigma_n(a)\sigma_n(x))=P(\sigma_1(ax),\ldots,\sigma_n(ax))=0
\end{equation*}
for all $x \in G$.  Taking a linear combination of $P$ and $Q$, we can find a polynomial with fewer terms than $P$ that vanishes on $\sigma_1x,\ldots,\sigma_nx$, giving a contradiction.
\end{proof}
{\it Proof of Theorem \ref{Hyper2}.}
Let $[M:k]=n$ and let $\sigma_1,\ldots,\sigma_n$ be the embeddings of $M$ into $\mathbb{C}$ fixing $k$.  Let $\alpha_0,\ldots,\alpha_{n-1}$ be a basis for $M$ over $k$.  Let $R$ be as in (\ref{Req}).

The only if direction has already been proven in the first half of our proof of Theorem \ref{Hyper}.  So suppose that there exists a nonzero homogeneous polynomial vanishing on $R$.  We need to show that $M$ contains a CM subfield over $k$.  By Lemma \ref{HyperL}, there exist nonidentical sequences of nonnegative integers $a_1,\ldots,a_n$ and $b_1,\ldots,b_n$ with $\sum_{i=1}^na_i=\sum_{i=1}^nb_i$ such that
\begin{equation}
\label{HyperE}
\prod_{i=1}^n \sigma_i(u)^{a_i}=\prod_{i=1}^n \sigma_i(u)^{b_i},  \qquad \forall u\in \co_M^*.
\end{equation}
By canceling terms, we can clearly assume that either $a_i=0$ or $b_i=0$ for $i=1,\ldots,n$.  Let $T$ be the set of $\sigma_i$'s such that $a_i\not=0$ and let $T'$ be the set of $\sigma_i$'s such that $b_i\not=0$.  By our assumption, $T$ and $T'$ are disjoint.  By composing both sides of (\ref{HyperE}) with some $\sigma_j$ we can assume that the identity embedding, $id$, is in $T$ (having fixed an identification of $M\subset \mathbb{C}$).  Let $\tau$ denote complex conjugation.  Let $\sigma_i \in T$.  We claim that $\sigma_j=\tau \sigma_i$ for some $\sigma_j\in T'$ and that $a_i=b_j$.  By the Dirichlet unit theorem, we can find a unit $u\in \co_M^*$ such that $|\sigma_i(u)|$ is very large and $|\sigma_l(u)|$ is very small and approximately the same size for all $\sigma_l\not=\sigma_i,\tau \sigma_i$.  Using that $\sum_{i=1}^na_i=\sum_{i=1}^nb_i$, this would clearly contradict (\ref{HyperE}) unless $\tau \sigma_i\in T'$  and $a_i=b_j$, where $\sigma_j=\tau \sigma_i$.  Applying the same argument to $T'$, we find that if $\sigma \in T'$ then $\tau \sigma \in T$.  Therefore $T'=\{\tau\sigma:\sigma \in T\}$.  In particular, $\tau \in T'$ and so $k$ must be real.  Since $T$ and $T'$ are disjoint, $T$ must consist only of complex embeddings.  

Let $\tilde{M}$ denote the Galois closure of $M$ over $k$.  Let $G=\Gal(\tilde{M}/k)$ and $H=\Gal(\tilde{M}/M)$.  Lift each $\sigma_i$ to an element $\tilde{\sigma}_i\in G$ such that  $\tilde{\sigma}_i|_M=\sigma_i$.  Let $\tilde{T}=TH=\{\tilde{\sigma}_ih:h\in H,i=1,\ldots,n\}$.  Similarly, let $\tilde{T}'=T'H$.  These definitions clearly don't depend on the liftings $\tilde{\sigma_i}$.  Of course, we still have $id_{\tilde{M}}\in\tilde{T}$, $\tau \in \tilde{T}'$, and $\tilde{T}$ and $\tilde{T}'$ are disjoint.  Let $\Sigma_{\tilde{M}}$ be the embeddings of $\tilde{M}$ into $\mathbb{C}$ (not necessarily fixing $k$).  Let $\phi \in \Sigma_{\tilde{M}}$.  Conjugating (\ref{HyperE}) by $\phi$ , we obtain
\begin{equation*}
\prod_{\sigma_i\in T} [\phi\sigma_i\phi^{-1}(u)]^{a_i}=\prod_{\sigma_i\in T'} [\phi\sigma_i\phi^{-1}(u)]^{b_i}=\prod_{\sigma_i\in T} [\phi\tau\sigma_{i}\phi^{-1}(u)]^{a_i}, \forall u\in \co_{\phi(M)}^*,
\end{equation*}
where the second equality follows from our earlier observations.  Note that each $\phi\sigma_i\phi^{-1}$ and $\phi\tau\sigma_{i}\phi^{-1}$ is an embedding of $\phi(M)$ into  $\mathbb{C}$ over $\phi(k)$.  Therefore, applying our previous reasoning to $\phi(M)$ and $\phi(k)$, we find that $\phi(k)$ is real (so $k$ is totally real) and that if $\sigma_i\in T$, then
\begin{equation}
\label{psp}
\phi\sigma_i\phi^{-1}=\tau\phi\tau\sigma_{j}\phi^{-1}
\end{equation}
on $\phi(M)$ for some $\sigma_j\in T$.  Since $\phi(\tilde{M})$ is Galois over $\phi(k)$ and $\phi(k)$ is real, $\tau \phi(\tilde{M})=\phi(\tilde{M})$.  It then makes sense to apply $\tau \phi^{-1} \tau$ on the left of each side of (\ref{psp}) to obtain $\tau \phi^{-1} \tau \phi \tilde{\sigma}_i\in \tilde{\sigma}_{j}H\subset \tilde{T}$.  So we see that 
\begin{equation*}
\tau \phi^{-1} \tau \phi\tilde{T}=\tilde{T}, \quad \forall \phi \in \Sigma_{\tilde{M}}.
\end{equation*}
Let $N=\langle \tau \phi^{-1} \tau \phi:\phi \in \Sigma_{\tilde{M}}\rangle$ be the subgroup of $G$ generated by the $\tau \phi^{-1} \tau \phi$'s.  Since $H$ is in $\tilde{T}$, we have in particular that $NH\subset \tilde{T}$.  Let $N'=\langle \tau\rangle N$.
\begin{lemma}
$N$  and $N'$ are normal subgroups of $G$.
\end{lemma}
\begin{proof}
Let $g\in G$ and $\phi \in \Sigma_{\tilde{M}}$.  By the definition of $N$ we see that 
\begin{equation*}
\tau(g^{-1}\tau)^{-1} \tau g^{-1}\tau=g\tau g^{-1} \tau\in N \text{ and }\tau (\phi g^{-1})^{-1} \tau \phi g^{-1}\in N.
\end{equation*}
Multiplying these two elements gives $g(\tau \phi^{-1} \tau \phi) g^{-1}\in N$ and therefore $N$ is a normal subgroup of $G$.  This implies $N'$ is actually a group, and as it is generated by $N$ and elements of the form $\phi \tau \phi^{-1}$, it is clearly a normal subgroup of $G$.
\end{proof}
Therefore $NH\subset \tilde{T}$ and $N'H$ are subgroups of $G$.  Let $L=\tilde{M}^{NH}$ be the fixed field of $NH$ and $L'=\tilde{M}^{N'H}$.  Since $M$ is the fixed field of $H$, we get inclusions $k\subset L'\subset L\subset M$.
\begin{lemma}
\label{CM}
$L$ is a CM field and $L'$ is its maximal real subfield.
\end{lemma}
\begin{proof}
Showing that $L$ is totally imaginary is equivalent to showing that for all $\phi \in \Sigma_{\tilde{M}}$, $\phi^{-1}\tau\phi \notin NH$.  If $\phi^{-1}\tau\phi \in NH$, then $\tau \in NH\subset \tilde{T}$, but since $\tau \in \tilde{T'}$, and $\tilde{T}$ and $\tilde{T'}$ are disjoint, we would have a contradiction.  Therefore $L$ is totally imaginary.  We now show that $L'$ is totally real.  This is equivalent to showing that $\phi^{-1}\tau\phi \in N'H, \forall \phi \in \Sigma_{\tilde{M}}$, which is trivial from the definition of $N'$.  Since we clearly have $[N'H:NH]=2$, we see that $L$ is a quadratic extension of $L'$.  Therefore $L$ is a CM field and $L'$ is its maximal real subfield.
\end{proof}
So we see that if $\co_M^*$ is not a dense subset of $M$ over $k$ then $M$ contains a CM subfield over $k$, and so the proofs of Theorems \ref{Hyper} and \ref{Hyper2} are complete.
\end{proof}
In fact, the field $L$ in Lemma \ref{CM} is the maximal CM subfield of $M$ over $k$.
\begin{lemma}
Let $M$ be a finite extension of a number field $k$.  Suppose that $M$ contains a CM subfield over $k$.  Then there exists a (unique) maximal CM subfield $L$ of $M$ over $k$, i.e.,  for any CM subfield $K$ of $M$ over $k$, $K\subset L$.  
\end{lemma}
\begin{proof}
Let $\tilde{M}$ be the Galois closure of $M$ over $k$.  Let $G=\Gal(\tilde{M}/k)$ and let $H=\Gal(\tilde{M}/M)$.  Let $\Sigma_M$ be the embeddings of $M$ into  $\mathbb{C}$.  Let $K$ be a CM subfield of $M$ over $k$ with maximal real subfield $K'$.  Let $\phi \in \Sigma_{\tilde{M}}$.  Let $\tau$ denote complex conjugation.  Then $\phi^{-1}\tau\phi$ gives an automorphism of $K$ over $K'$ since $K$ is a CM field.  Since $K$ is totally imaginary, this automorphism cannot be the identity on $K$.  Therefore it is complex conjugation, and so $\tau\phi^{-1}\tau\phi$ fixes $K$, that is $\tau\phi^{-1}\tau\phi \in \Gal(M/K)$.  Let $N\subset G$ be the group generated by the $\tau\phi^{-1}\tau\phi$'s.  Since $H\subset \Gal(M/K)$, we have $NH\subset \Gal(M/K)$.  Since $K$ is complex, $\tau\notin NH$.  But then the proof of Lemma \ref{CM} shows that the fixed field of $NH$, $L$, is a CM subfield of $M$ over $k$, and by Galois theory $K\subset L$.  So $L$ is the maximal CM subfield of $M$ over $k$.
\end{proof}
\section{Non-archimedean Places}
We now consider the general case, where $S$ may contain non-archimedean places.  Of course, when Theorem \ref{Hyper} tells us that the $S_\infty$-integral points on $\mathbb{P}^n\backslash Z$ are dense, it is trivial that the $S$-integral points on $\mathbb{P}^n\backslash Z$ are dense for any $S$ (containing $S_\infty$).  Furthermore, the proof assuming that (a) is satisfied works for arbitrary $S$, and condition (b) doesn't occur if $S$ contains non-archimedean places.  The real difficulty arises when condition (c) of Theorem~\ref{Hyper} occurs and $S$ is larger than $S_\infty$.

Assuming that neither (a) or (b) of Theorem \ref{Hyper} holds and that (c) is satisfied, we easily reduce, as before, to the case where $Z\subset \mathbb{P}^{m-1}$ is irreducible over $k$ defined by $N^M_k(x_0\alpha_0+\cdots+x_{m-1}\alpha_{m-1})=0$, where $M$ contains $L$, the maximal CM subfield of $M$ over $k$.  Using Lemma \ref{units}, the problem is equivalent to determining if $\co_{M,S_M}^*$ is a dense subset of $M$ over $k$, where $S_M$ is the set of places of $M$ lying over places of $S$.  Thus, we are in a position to apply Lemma \ref{HyperL}.  Paying careful attention to the proof of Theorem \ref{Hyper}, we see that if (\ref{eq}) holds for all $x\in \co_{M,S_M}^*$, then the identity must be of the form
\begin{equation*}
\prod_{i=1}^l \sigma_i N^M_L(x)^{a_i}=\prod_{i=1}^l  \tau \sigma_i N^M_L(x)^{a_i}, \quad \forall x\in \co_{M,S_M}^*,
\end{equation*}
where $\tau$ denotes complex conjugation and $\sigma_1,\ldots,\sigma_l$ are the embeddings of $L$ into $\mathbb{C}$ fixing $k$.  Since $N^M_L(\co_{M,S_M}^*)$ is a finite index subgroup of $\co_{L,S_L}^*$, we can reduce to the problem of determining whether there is a nontrivial identity
\begin{equation}
\label{eq2}
\prod_{i=1}^l (\sigma_ix)^{a_i}=\prod_{i=1}^l   (\tau\sigma_ix)^{a_i}
\end{equation}
for all $x\in \co_{L,S_L}^*$.  Without loss of generality, by raising both sides of (\ref{eq2}) to an appropriate power, we can assume that $a_i$ is divisible by $[\co_{L,S_L}^*:\co_{L',S_{L'}}^*]$ for all $i$, where $L'$ is the maximal real subfield of $L$.  In that case, (\ref{eq2}) is true for all $x\in \co_L^*$ and any appropriate choice of the $a_i$.  So we can essentially reduce to studying $\co_{L,S_L}^*/\co_{L}^*$.  Assume now that $L/k$ is Galois with Galois group $G=\Gal(L/k)$.  If one can compute a minimal set of generators for the free abelian group $\co_{L,S_L}^*/\co_{L}^*$ and the action of $G$ on it in terms of those generators, then determining the existence of a solution to (\ref{eq2}) becomes elementary linear algebra. So, at least in the case the appropriate field extensions are Galois, the problem of determining whether there exists a Zariski-dense set of $S$-integral points on a complement of hyperplanes is reduced to being able to do certain computations with the non-archimedean part of the $S$-unit group in particular CM fields.

Of course, the action of $G$ on $\co_{L,S_L}^*/\co_{L}^*$ is closely related to how the non-archimedean places in $S$ split in $L$.  For instance (still assuming $L/k$ Galois), if some place of $S$ splits completely in $L$, then $\co_{L,S_L}^*$ is a dense subset of $L$ over $k$ (see the proof of Lemma \ref{HyperL}).  On the other hand, if no place of $S_{L'}$ splits in $L$, then  $\co_{L,S_L}^*$ is not a dense subset of $L$ over $k$.  More generally, let $D$ be the set of decomposition fields of the non-archimedean places of $S_L$.  Then it is easily shown that if ${L'}^*\prod_{F\in D}F^*$ is not a dense subset of $L$ over $k$, then $\co_{L,S_L}^*$ is not a dense subset of $L$ over $k$.  This leads to the following natural question.
\begin{question}
\label{qu}
Let $L$ be a finite extension of a number field $k$.  Let $\mathcal{F}$ be a set of subfields of $L$ over $k$.  Can one simply characterize when $\prod_{F\in \mathcal{F}}F^*$ is a dense subset of $L$ over $k$?
\end{question}
While this question does not seem to have been studied before, the related problems of determining when $\prod_{F\in \mathcal{F}}F^*=L^*$ and, more generally, determining the group structure of $L^*/\prod_{F\in \mathcal{F}}F^*$ have been studied in \cite{Wie},\cite{Wie3}, and \cite{Wie2}.  It would be interesting to connect this work to Question \ref{qu}.
\bibliography{integral}
\end{document}